\documentstyle[amssymb]{amsart}
\input mssymb 

\reversemarginpar

\reversemarginpar

\reversemarginpar

\def\grayrule{\special{ps:: currentgray 0.9 setgray}
\hbox to 0pt{\hss
\vrule height 0.3 in depth 0.2 in width 0.48 in \hss}
\special{ps:: setgray}}
\def\sqquare#1#2{\square\hbox to 0pt{$#1#2$\hss}}

\newcommand{\ZFCa}{\text{\normalshape\sf ZFC}}
\newcommand{\CH}{\text{\normalshape\sf CH}}

\newcommand{\lbv}{[\![}
\newcommand{\rbv}{]\!]}
\newcommand{\reals}{\Bbb R}
\newcommand{\rationals}{\Bbb Q}
\newcommand{\rest}{{\mathord{\restriction}}}
\newcommand{\add}{\text{\normalshape\sf{add}}}
\newcommand{\cov}{\text{\normalshape\sf{cov}}}
\newcommand{\unif}{\text{\normalshape\sf{non}}}
\newcommand{\Borel}{\text{\normalshape\sf{Borel}}}
\newcommand{\cof}{\text{\normalshape\sf{cof}}}

\newcommand{\ran}{{\text{\normalshape\sf {range}}}}

\newcommand{\QED}{\hspace{0.1in} \Box \vspace{0.1in}}
\newcommand{\forces}{\mathrel{\|}\joinrel\mathrel{-}}
\newcommand{\F}{{\cal F}}
\newcommand{\D}{{\cal D}}
\newcommand{\N}{{\cal N}}
\newcommand{\M}{{\cal M}}
\newcommand{\V}{{\bold V}}
\newcommand{\<}{\langle}
\renewcommand{\>}{\rangle}

\newcommand{\thinks}{\models}

\newtheorem{theorem}{Theorem}[section]
\newtheorem{lemma}[theorem]{Lemma}

\theoremstyle{definition}
\newtheorem{definition}[theorem]{Definition}

\newcommand{\lesdot}{\mathrel{\mathord{<}\!\!\raise 
0.8 pt\hbox{$\scriptstyle\circ$}}}
\newcommand{\Proof}{{\sc Proof} \hspace{0.2in}}

\newcommand{\lft}[2]{\mathopen\ifcase#1{}\oo\or
                        \big#2\or\Big#2\else\oo\fi} 
\newcommand{\rgt}[2]{\mathclose\ifcase#1{}\oo\or
                        \big#2\or\Big#2\else\oo\fi} 

\def\setup#1 {\gdef\firstparameter{#1}}
\def\ldiagram#1 #2 #3 #4 #5 #6 #7 #8 #9 {\vbox{\hsize=5.5 cm
\vsize=4 cm\relax
$$\begin{array}{ccccccc}
\firstparameter&
\rightarrow &
#1 &
\rightarrow &
#2 &
\rightarrow &
#3\\
\ &
\ &
\uparrow &
\ &
\uparrow &
\ &
\ \\
\smash{\bigg\uparrow} &
\ &
#4 &
\ &
#5 &
\ &
\smash{\bigg\uparrow} \\
\ &
\ &
\uparrow &
\ &
\uparrow &
\ &
\ \\
#6 &
\rightarrow &
#7 &
\rightarrow &
#8 &
\rightarrow &
#9
\end{array}$$}}
\def\cdiagram #1 #2 #3 #4 #5 #6 #7 #8 #9 {
\par
\hbox to \hsize{\hfill 
\ldiagram  #1 #2 #3 #4 #5 #6 #7 #8 #9  
\hfill}}

\def\ccdiagram{%
\hbox to \hsize{\hfill
\lldiagram
\hfill}}
\def\onesetup#1 #2 #3 #4 #5 {%
\gdef\pone{#1}
\gdef\ptwo{#2}
\gdef\pthree{#3}
\gdef\pfour{#4}
\gdef\pfive{#5}}
\def\twosetup#1 #2 #3 #4 #5 {%
\gdef\psix{#1}
\gdef\pseven{#2}
\gdef\peight{#3}
\gdef\pnine{#4}
\gdef\pten{#5}}

\def\lldiagram{{\vbox{\hsize=11 cm
\vsize=4 cm\relax
$$\begin{array}{ccccccc}
\pone&
\rightarrow &
\ptwo &
\rightarrow &
\pthree &
\rightarrow &
\pfour\\
\ &
\ &
\uparrow &
\ &
\uparrow &
\ &
\ \\
\smash{\bigg\uparrow} &
\ &
\pfive &
\rightarrow  &
\psix &
\ &
\smash{\bigg\uparrow} \\
\ &
\ &
\uparrow &
\ &
\uparrow &
\ &
\ \\
\pseven &
\rightarrow &
\peight &
\rightarrow &
\pnine &
\rightarrow &
\pten
\end{array}$$}}}

\renewcommand{\t}{{\cal T}}

\newcommand{\p}{{\cal P}}
\newcommand{\incomp}{\perp}
\newcommand{\compatible}{\not\perp}
\begin{document}
\title{Adding one random real}
\author{Tomek Bartoszy\'{n}ski}
\thanks{First author partially supported by  SBOE grant 92--096}
\address{Department of Mathematics\\
Boise State University\\
Boise, Idaho 83725 U.S.A.}
\email{tomek@@math.idbsu.edu}
\author{Andrzej Ros{\l}anowski}
\thanks{Second author partially supported by KBN grant 1065/P3/93/04}
\address{Department of Mathematics\\
Hebrew University\\
Jerusalem, Israel, and\\
Mathematical Institute\\
Wroclaw University\\
Wroclaw, Poland}
\email{roslanow@@sunrise.huji.ac.il}
\author{Saharon Shelah}
\thanks{Third author partially supported by Basic Research Fund,
Israel Academy of Sciences, publication 490}
\address{Department of Mathematics\\
Hebrew University\\
Jerusalem, Israel}
\email{shelah@@sunrise.huji.ac.il}
\maketitle
\begin{abstract}
We study the cardinal invariants of measure and category after adding
one random real. In particular, we show that the number of measure
zero subsets of the plane which are necessary to cover graphs of
all continuous functions maybe large while the covering for measure
is small.
\end{abstract}
\section{Introduction}
Let ${\cal J}$ be an ideal of
subsets of the real line (where real line means $\reals$, $2^\omega$
or $[0,1]$).
Define the following  cardinal invariants:
\begin{enumerate}

\item $\add({\cal J})=\min\{|{\cal A}|: 
{\cal A} \subseteq {\cal J} \ \&\ \bigcup {\cal A} \not \in
{\cal J}\}$,
\item $\cov({\cal J})=\min\{|{\cal A}|: 
{\cal A} \subseteq {\cal J} \ \&\ \bigcup {\cal A} =\reals\}$,
\item $\unif({\cal J})=\min\{|X|: X \subseteq \reals \ \& \ X \not\in
  {\cal J}\}$,
\item $\cof({\cal J})=\min\{|{\cal A}|: {\cal A} \subseteq {\cal J}
  \&\ \forall A \in {\cal J} \ \exists B \in {\cal A} \ A \subseteq B\}$.
\end{enumerate}

Let ${\cal M}$ and ${\cal N}$ be the ideals of meager and of measure zero
subsets of the real line respectively.
Finally let
${\frak b}$ be the size of the smallest
unbounded family in $\omega^{\omega}$ and
${\frak d}$ the size of the smallest
dominating family in $\omega^{\omega}$ .

The relationship between these cardinals is
described in the following diagram, where arrows
  means $\leq $:

\onesetup  \cov(\N) \unif(\M) \cof(\M) \cof(\N) {\goth b}

\twosetup {\goth d} \add(\N) \add(\M) \cov(\M) \unif(\N)

\ccdiagram

In addition $\add({\cal M}) = \min\{{\frak b}, \cov({\cal M})\}$ and
$\cof({\cal M}) = \max\{{\frak d}, \unif({\cal M})\}$.

The proofs of those inequalities can be found in \cite{Bar84Add}, \cite{fre:cichon}
and \cite{Mil81Som}. 
In this paper we show that except for $\cov(\N)$ and $\unif(\N)$
values of these invariants do not change when one random real is added.
Let $\bold B $ be the measure algebra adding one random real.
\begin{theorem}[Pawlikowski, Krawczyk]\label{pawlikowski}
The following holds in $\V^{\bold B}$:
\begin{enumerate}
\item $\add(\N)=\add(\N)^{\bold V}$ and 
$\cof(\N)=\cof(\N)^{\V}$,
\item $\cov(\N) \geq \cov(\N)^{\V}$ and 
$\unif(\N)\leq\unif(\N)^{\V}$,
\item $\cov(\N) \geq {\frak b}$ and $\unif(\N) \leq {\frak d}$,
\item  ${\frak b} = {\frak b}^{\V}$ and ${\frak d} = {\frak d}^{\V}$,
\item $\cov(\M) \geq \cov(\M)^{\V}$ and
$\unif(\M) \leq \unif(\M)^{\V}$,
\item $\add(\M) \geq \add(\M)^{\V}$ and 
$\cof(\M) \leq \cof(\M)^{\V}$.
\end{enumerate}
\end{theorem}
\Proof (1), (2) and (4) is folklore (see \cite{Paw86Sol}).
(3) is due to Krawczyk (see \cite{Paw86Sol} or \cite{Kra83}), (5) is
due to Pawlikowski 
(\cite{Paw86Sol}) and (6) follows from (5), (4) and the remarks above.~$\QED$

For a set $H \subseteq \reals\times\reals$ and $x,y \in \reals$ let
$(H)_x = \{y : \<x,y\> \in H$ and let $(H)^y=\{x : \<x,y\> \in H\}$.

We will use the following classical lemma:
\begin{lemma}\label{represent}
  Suppose that $r$ is a random real over $\V$. Then
  \begin{enumerate}
    \item for every $x \in \V[r] \cap \reals$ there exists a Borel
      function $f : \reals \longrightarrow \reals$ such that $x =
      f(r)$,
    \item for every Borel measure zero set $F \in \V[r]$ there exists
      a Borel measure zero set $H \subseteq \reals\times\reals$, $H
      \in \V$ such that $F = (H)_r$.~$\QED$  
  \end{enumerate}
\end{lemma}
We will need the following characterization of $\cov(\M)$ and
$\unif(\M)$.

Let $$\cal S = \{S \in ([\omega]^{<\omega})^\omega : \forall n \ |S(n)|
\leq (n+1)^2\}.$$

\begin{theorem}[\cite{Bar87Com}]\label{cohen}
  The following conditions are equivalent:
  \begin{enumerate}
  \item $\cov(\M)\geq\kappa$,
  \item for every family $F \subseteq \omega^\omega$ of size $<\kappa$
    there exists $g \in \omega^\omega$ such that 
$$\forall f \in F \ \exists^\infty n \ f(n)=g(n).$$
\item for every family $F \subseteq \omega^\omega$ of size $<\kappa$
    there exists $S \in \cal S$ such that 
$$\forall f \in F \ \exists^\infty n \ f(n) \in S(n).$$
  \end{enumerate}
Similarly,
\begin{enumerate}
\item $\unif(\M)\geq\kappa$,
\item for every family $F \subseteq \omega^\omega$ of size $<\kappa$
    there exists $g \in \omega^\omega$ such that 
$$\forall f \in F \ \forall^\infty n \ f(n) \neq g(n).$$
\item for every family $F \subseteq \cal S$ of size $<\kappa$
    there exists $f  \in \omega^\omega $ such that 
$$\forall S \in F \ \forall^\infty n \ f(n) \not \in S(n).~\QED $$
\end{enumerate}
\end{theorem}

\section{Cohen reals}
In this section we will show that invariants $\cov(\M)$ and
$\unif(\M)$ do not change when random reals are added.
\begin{theorem}
  The following holds in $\V^{\bold B}$:
\begin{enumerate}
\item $\cov(\M) = \cov(\M)^{\V}$ and
$\unif(\M) = \unif(\M)^{\V}$,
\item $\add(\M) = \add(\M)^{\V}$ and 
$\cof(\M) = \cof(\M)^{\V}$.
\end{enumerate}
\end{theorem}
\Proof (1) By \ref{pawlikowski}, it is enough to show that in   
$\V^{\bold B}$, $\cov(\M) \leq \cov(\M)^{\V}$ and
$\unif(\M) \geq \unif(\M)^{\V}$. 

By \ref{cohen}, there exists a family $F \subseteq \omega^\omega$ of
size $\cov(\M)^{\V}$ such that 
$$\forall S \in {\cal S} \ \exists f \in F \ \forall^\infty n \ f(n)
\not \in S(n).$$

By \ref{cohen}, to finish the proof it is enough to show that
$$\V^{\bold B} \thinks \forall g  \in \omega^\omega  \ \exists f \in F
\ \forall^\infty n \ f(n) \neq g(n).$$

Let $\dot{g}$ be a $\bold B$-name for an element of $\omega^\omega$.
Define for $n \in \omega$,
$$S(n) = \left\{k \in \omega : \mu\left(\lbv \dot{g}(n)=k \rbv_{\bold
B}\right)>\frac{1}{(n+1)^2} \right\},$$
where $\mu$ is the Lebesgue measure.
It is clear that $|S(n)| <(n+1)^2$ for all $n$. Therefore there exists $f
\in F$ and $N \in \omega$ such that $f(n) \not \in S(n)$ for all $n\geq N$.
We claim that 
$$\forces_{\bold B} \forall^\infty n \ f(n) \neq \dot{g}(n).$$
Let $p \in \bold B$. Find $n>N$ such that $\sum_{k=n}^\infty k^{-2} <
\mu(p)$. Then
$$q = p - \bigcup_{k=n}^\infty \lbv \dot{g}(k) = f(k)\rbv_{\bold B}
>0$$
and 
$$q \forces_{\bold B} \forall k>n \ f(k) \neq \dot{g}(k).$$

\vspace{0.1in}
 
To show that 
$\unif(\M) \geq \unif(\M)^{\V}$ holds in $\V^{\bold B}$, we
``dualize'' the above argument. 

Suppose that $F \subseteq \omega^\omega$ is a family of size
$\unif(\M)$ in $\V^{\bold B}$ such that 
$$\V^{\bold B} \thinks \forall g \in \omega^\omega \ \exists f \in F \
\exists^\infty n \ f(n)=g(n).$$
Let $\dot{F}=\{\dot{f} : f \in F\}$ be a set of $\bold B$-names for
elements of $F$. Without loss of generality we can assume that
$\dot{F} \in \V$. For $\dot{f} \in \dot{F}$ let $S_f \in \cal S$ be
defined as 
$$S_f(n) = \left\{k \in \omega : \mu\left(\lbv \dot{f}(n)=k \rbv_{\bold
B}\right)>\frac{1}{(n+1)^2} \right\}.$$
As before we show that 
$$\forall g \in \omega^\omega \ \exists \dot{f} \in \dot{F} \
\exists^\infty n \ g(n) \in S_f(n),$$
which by \ref{pawlikowski}, finishes the proof.

To show the second part use \ref{pawlikowski} and the fact that
$\add({\cal M}) = \min\{{\frak b}, \cov({\cal M})\}$ and
$\cof({\cal M}) = \max\{{\frak d}, \unif({\cal M})\}$.~$\QED$

Recall that a set $X \subseteq \reals$ has strong measure zero if for
every sequence of positive reals  $\<\varepsilon_n: n \in \omega\>$
there exists a sequence of intervals $\<I_n :n \in \omega\>$ such that
the length of $I_n$ is $\leq \varepsilon_n$ and $X \subseteq
\bigcup_{m} I_m$. Note that, equivalently we can request that
$X \subseteq
\bigcap_{n \in \omega} \bigcup_{m>n} I_m$.
\begin{theorem}\label{strmeazer}
  Suppose that $X \subseteq \reals$ and $X \in \V$. Then  $X$ has
  strong measure zero in $\V$ iff $X$ has
  strong measure zero in $\V^{\bold B}$.
\end{theorem}
\Proof
It is easy to see that for every sequence $\<\varepsilon_n: n \in
\omega\> \in \V^{\bold B}$ there exists a sequence $\<\delta_n : n \in
\omega\> \in \V$ such that $\delta_n \leq \varepsilon_n$ for all $n$.
Therefore, if $X$ has
  strong measure zero in $\V$ then $X$ has
  strong measure zero in $\V^{\bold B}$.

Suppose that $X$ does not have strong measure zero in $\V$ and let
$\<\varepsilon_n : n \in \omega\>$ be a sequence of positive reals
witnessing that.  Suppose that $X$ has strong measure zero in $\V^{\bold
B}$. Let $\<\delta_n : n \in \omega\>$ be a decreasing sequence of positive
reals such that $\delta_{n} <\varepsilon_k$ for all $k\leq n^3$.  Let
$\<\delta_n':n\in\omega\>$ be a decreasing sequence of positive rationals
such that $\delta_{2k}'=\delta_{2k+1}'$ and $\delta_n'<\delta_n$. By the
assumption we can find a sequence of intervals $\<I_n : n \in
\omega\> \in \V^{\bold B}$ such that $X \subseteq 
\bigcap_{n \in \omega} \bigcup_{m>n} I_m$ and the length of $I_m$ is less
than $\delta_{2m}'$.
Let $\<I(n,k):k \in \omega\>$ be a partition of $\reals$ into rational intervals
of the length $\delta_n'$. Each interval $I_m$ is covered by $I(2m,k)\cup
I(2m,k+1)$ for some $k=k(m)$. Let $\<\dot{I}_n:n\in\omega\>$ be a ${\bold
B}$-name for the sequence $\<I(2m,k(m)), I(2m,k(m)+1):m\in\omega\>$ (i.e.
$\dot{I}_{2m}$ is a name for $I(2m,k(m))$ and $\dot{I}_{2m+1}$ is that for
$I(2m, k(m)+1)$). Thus
$$\forces_{\bold B}\text{``the length of $\dot{I}_n$ is $\delta_n'$ \&
$X\subseteq\bigcap_{n\in\omega}\bigcup_{m>n}\dot{I}_m$''.}$$
Now define for $n=2m+i$ ($i=0,1$):
$${\cal A}_n = \left\{I(2m,k) : \mu\left(\lbv\dot{I}_n=I(2m,k) \rbv_{\bold 
B}\right)>\frac{1}{(n+1)^2} \right\}.$$
Note that $|{\cal A}_n|<(n+1)^2$ (some ${\cal A}_n$'s may be empty). By the
choice of the sequence $\<\delta_n : n \in \omega\>$ if we order
lexicographically the intervals in $\bigcup_{n \in \omega} {\cal A}_n$ in a
sequence $\<J_n : n \in \omega\>$, then the length of $J_n$ will be $\leq
\varepsilon_n$. Let $x \in X$ be such that $x \not \in
\bigcup_{n \in \omega} J_n$. Note that then for each $n\in\omega$
$$\mu(\lbv x\in\dot{I}_n\rbv)\leq\frac{1}{(n+1)^2}$$
Let $p \in \bold B$. Find $n$ such that $\sum_{k=n}^\infty k^{-2} <
\mu(p)$. Then
$$q = p - \bigcup_{k=n}^\infty \lbv x \in \dot{I}_n\rbv_{\bold B}
>0$$
and 
$$q \forces_{\bold B} \forall k>n \ x \not\in \dot{I}_n.~\QED$$

The proof of \ref{strmeazer} seems to suggest that a filter $\F\in\V$ on
$\omega$ which cannot be extended to a rapid filter in $\V$ cannot be
extended to a rapid filter in $\V^{\bold B}$. However, this is not the
case. First, let us recall that a non-principal filter $\F$ on $\omega$ is
called rapid if for every increasing function $f\in\omega^\omega$ there
exists $X\in\F$ such that $|X\cap f(n)|\leq n$ for all $n$. 
\begin{theorem}
Suppose that $\D$ is a rapid filter on $\omega$. Then there
exists a filter $\F$ such that:
\begin{enumerate}
\item $\F$ cannot be extended to a rapid filter in $\V$,
\item $\V^{\bold B}\models$``$\F$ can be extended to a rapid filter''.
\end{enumerate}
\end{theorem}
\Proof Let $\F$ be the family of all subsets $A$ of $\omega$ such that for
some set $X\in\D$ the sequence 
$$\frac{|A\cap [n^2,(n+1)^2)|}{2n+1} \stackrel{n\in
  X}{\longrightarrow} 1.$$
It should be clear that if $A\subseteq B$, $A\in \F$ then
$B\in \F$ and the same set $X\in\D$ witnesses it. Moreover if $A,B\in\F$ is
witnessed by $X_A,X_B\in \D$ then the intersection $X_A\cap X_B\in\D$
witnesses that $A\cap B\in \F$. Consequenlty $\F$ is a non-principal
filter on $\omega$. We claim that $\F$ cannot be extended to a rapid
filter. 
Suppose that a set $A\subseteq\omega$ is such that $|A\cap n^3|\leq n$
for $n \in \omega$. Then for each $m\in\omega$ we have
$$\frac{|A\cap [m^2,(m+1)^2)|}{2m+1}\leq \frac{[(m+1)^{2/3}]+1}{2m+1}\leq
\frac{(m+1)^{2/3}+1}{2m+1}$$
and hence 
$$\lim_{m\rightarrow\infty}\frac{|A\cap [m^2,(m+1)^2)|}{2m+1}=0.$$
Consequently the complement $\omega\setminus A$ of the set $A$ belongs to
$\F$ and $A$ cannot be in any filter extending $\F$. 

\vspace{0.2in}

To prove the assertion $(2)$ we work with the measure algebra on the
space $\prod\limits_{n\in\omega}[n^2,(n+1)^2)$ equipped with
the natural product measure $\mu$ (we use the same symbol as for the
Lebesgue measure since this measure corresponds to Lebesgue measure
under cannonical mapping of underlying space onto the interval
$[0,1]$).
Suppose that
$r\in\prod\limits_{n\in\omega} 
[n^2,(n+1)^2)$ is a random real over $\V$ and work in $\V[r]$. 
First note
that for a set $A \in \F$ and $X \in \cal D$,
\begin{multline*}
\mu\left(\left\{x \in  \prod\limits_{n\in\omega}[n^2,(n+1)^2):
\forall^\infty m\in X \ x(m)\notin A\right\}\right) \leq \\
\sum_{n=0}^\infty
\prod^\infty\begin{Sb}m=n\\
m \in X
\end{Sb}
 \left(1-\frac{|A\cap [m^2, (m+1)^2)|}{2m+1}\right)=0.
\end{multline*}

In particular, since $r$ is a random real,
$$\forall X \in \cal D \ \forall A \in \F \ A \cap \ran(r\rest X) \neq
\emptyset.$$

Consequently $\F\cup\{\ran(r\rest X): X \in \cal D\}$ generates a
filter $\F^\star$. We are going 
to show that it is a rapid filter. Suppose that $f\in\omega^\omega\cap
\V[r]$ is an increasing function. Since random real forcing is
$\omega^\omega$-bounding we can assume that $f \in \V$.  Since $\D$
was a rapid filter in $\V$ we find 
a set $X\in D$ such that $|X\cap f(n)|\leq n$ for $n \in \omega$. Look
at the set $A=\{r(n) : n \in X\}$. For
every $n\in\omega$ we have:
$$|A\cap f(n)|\leq  |X\cap f(n)|
\leq n.$$
The theorem is proved.~$\QED$

\section{Random reals}
Theorem \ref{pawlikowski} shows that in $ \V^{\bold B}$,  $\cov(\N) \geq
\max\{\cov(\N)^{\V}, {\frak b}^{\V}\}$. In this section we will show
that it is consistent that $\cov(\N)^{\V^{\bold B}} >
\max\{\cov(\N)^{\V}, {\frak b}^{\V}\}$. 

We will need the following notation:
\begin{definition}
  Let $\N_2$ be the ideal of measure zero subsets of
$\reals\times\reals$ and let $\Borel(\reals)$ be the collection of
all Borel mappings from $\reals $ into $\reals$. Define

\begin{multline*}
\cov^\star(\N)=\min\lft2\{|{\cal A}|:  {\cal A} \subseteq \N_2 \ \&\
\forall f \in \Borel(\reals) \ \forall B \in \Borel \setminus \N \ \exists
H \in {\cal A} \\
 \lft1\{x \in B: \<x,f(x)\> \in H\lft1\}\not\in \N\rgt2\}
\end{multline*}
 and 
\begin{multline*}
\unif^\star(\N)=\min\lft2\{|X| : X \subseteq \Borel(\reals) \ \& \ \forall
H \in \N_2 \ \forall B \in \Borel \setminus \N \ \exists f \in X \\
 \lft1\{x
\in B : \<x,f(x)\> \not \in H\rgt1\} \not\in \N\rgt2\}.
\end{multline*}
\end{definition}

As a consequence of \ref{represent}, we get:
\begin{lemma}
  $\cov^\star(\N)=\cov(\N)^{\V^{\bold B}}$ and
  $\unif^\star(\N)=\unif(\N)^{\V^{\bold B}}$.~$\QED$
\end{lemma}

The goal of this section is to show that the coefficient
$\cov^\star(\N)$ can be large while both $\frak b$ and $\cov(\N)$ are
small and that $\unif^\star(\N)$ can be small while both $\unif(\N)$
and $\frak d$ are large.

The key to our construction is the following theorem:
\begin{theorem}
  There exists a forcing notion ${\cal P}$, adding generically a
  continuous function $h_G:\reals \longrightarrow \reals$, such that  
  \begin{enumerate}
    \item ${\cal P}$ is $\sigma$-centered,
    \item $\forall f \in \omega^\omega \cap \V^{\cal P} \ \exists g
      \in \V \cap \omega^\omega \ \exists^\infty n \ f(n) \leq g(n)$,
    \item for every $H \in \N_2 \cap \V$, $\{x: \<x,h_G(x)\> \in H\}$
      has measure zero.
  \end{enumerate}
\end{theorem}
\Proof
Let $\t$ consists of all pairs $\<\varepsilon,\phi\>$ where $\varepsilon$ is a
rational number in $(0,1)$ and $\phi:2^{<\omega} \times 2^{<\omega}
\longrightarrow [0,1]$ is a 
function such that for $s,t \in 2^{<\omega}$:
\begin{enumerate}
\item $\phi(\emptyset,\emptyset)>0$,
\item $\phi(s,t)\leq 2^{-(|s|+|t|)}$,
\item $\phi(s^\frown 0,t)+\phi(s^\frown 1,t)
=\phi(s,t)=\phi(s,t^\frown 0)+\phi(s,t^\frown 1)$.
\end{enumerate}

We define the partial order ${\cal P}$. Conditions are pairs $p=\<h,u\>$ such
that 
\begin{enumerate}
\item $u\in [\t]^{<\omega}$
\item $h:2^{\leq m}\longrightarrow 2^{<\omega}$ for some $m=m(p)$,
\item if $s \subseteq t \in 2^{\leq m}$ then $h(s)\subseteq h(t)$,
\item if $\<\varepsilon,\phi\>\in u$ then
$$\sum_{s\in
2^m} 2^{|h(s)|} \phi\lft1(s,h(s)\rgt1) > \varepsilon.$$ 
\end{enumerate}

The order $\leq$ on ${\cal P}$ is the natural one:
$$(h,u)\geq (h',u')\iff h\supseteq h'\ \&\ u\supseteq u'.$$

\begin{lemma}
\label{density}
Suppose that $p=(h,u)\in{\cal P}$. Then there is $q=(h',u)\in{\cal P}$
such that $q\geq p$, $m(q)>m(p)$ and if $s \in 2^{m(q)}$ then $
\lft1|h'(s)\rgt1|> \lft2|h'\lft1(s\rest m(p)\rgt1)\rgt2|$.
\end{lemma}
\Proof What we have to do is to extend $h$. Note that if we put
$h'(s^\frown i)=h(s)$ (for $s \in 2^{m(p)}$) then $(h',u)$ is a
condition stronger than $p$.  So the only
problem is to extend the ``values'' of $h$.

Take $\delta>0$ such that for every $\<\varepsilon,\phi\>\in u$
$$\delta\cdot\sum_{s \in 2^{m(p)}}2^{1+|h(s)|}< \sum_{s\in
2^m(p)} 2^{|h(s)|} \phi\lft1(s,h(s)\rgt1)-\varepsilon.$$

\begin{lemma}
There are $m'>m(p)$ and $e:2^{m'}\longrightarrow
2$ such that for each $\<\varepsilon,\phi\>\in u$, $s\in 2^{m(p)}$:
$$\frac{1}{2}\phi\lft1(s,h(s)\rgt1)-\delta<\sum_{s\subseteq t
\in 2^{m'}} \phi\lft1(t,h(s)^\frown e(t)\rgt1).\leqno(\otimes)$$ 
\end{lemma}
\Proof Let $n=|u|$ and let $m'>m(p)$ be such that $2^{-m'}/\delta^2<1/n$.
Fix $s\in 2^{m(p)}$. We are going to find a function $e_s:\{t\in
2^{m'}:s\subseteq t\}\longrightarrow 2$ satisfying the condition
$(\otimes)$ for each $\<\varepsilon,\phi\>\in u$. Consider the space
$\Omega$ of all functions from $\{t\in 2^{m'}: s\subseteq t\}$ to $2$. The
space carries the natural (product) probability measure $P$. For
$\<\varepsilon,\phi\>\in u$ define a random variable
$Y_\phi:\Omega\longrightarrow [0,1]$ by 
$$Y_\phi(e)=\sum_{s\subseteq t\in 2^{m'}}\phi\lft1(t,h(s)^\frown e(t)\rgt1).$$
By the Tchebyshev inequality  we know that 
$$P\left(\left|Y_\phi-\int
Y_\phi\, d \Omega
\right|\geq\delta\right)\leq\frac{\bold D^2 Y_\phi}{\delta^2}.$$ 
If we put $X^t_\phi(e)=\phi\lft1(t,h(s)^\frown e(t)\rgt1)$ (for $t\in 2^{m'},
s\subseteq t$) then $X^t_\phi$'s are independent random variables on
$\Omega$ and $Y_\phi=\sum\limits_{s\subseteq t\in 2^{m'}} X^t_\phi$. Now,
\begin{multline*}
\bold D^2
Y_\phi=\int \left(Y_\phi-\int Y_\phi\, d \Omega\right)^2\, d\Omega
= 
\int\left(\sum_{s\subseteq t\in
2^{m'}}\left(X^t_\phi-\int X^t_\phi\, d \Omega\right)\right)^2\, d\Omega =\\
\sum_{s\subseteq t\in 
2^{m'}}\int
\left(X^t_\phi-\int X^t_\phi \, d\Omega\right)^2\, d\Omega
\end{multline*} 
(for the last equality we use the independence of $X^t_\phi$'s). Since
$\left|X^t_\phi -\int X^t_\phi\, d\Omega\right|\leq
2^{-(m'+|h(s)|+1)}$ we get  
$$\bold D^2 Y_\phi\leq 2^{m'-m(p)}\cdot 2^{-2m'-2|h(s)|-2}<2^{-m'}.$$
Hence
$$P\left(\left|Y_\phi-\int Y_\phi\, d\Omega \right|\geq\delta\right)<2^{-m'}/\delta^2<1/n$$
and therefore we can find $e_s\in\Omega$ such that for each
$(\varepsilon,\phi)\in u$ we have $\int  Y_\phi - \delta\, d\Omega \leq
Y_\phi(e_s)$. Since
$$\int  Y_\phi\, d\Omega =\sum\limits_{s\subseteq t\in
2^{m'}}\int  X^t_\phi\, d\Omega = \frac{1}{2}\sum\limits_{s\subseteq
t\in 2^{m'}} \phi(t,h(s))= \frac{1}{2}\phi(s,h(s))$$
we get that $e_s$ is as required.~$\QED$

Define $h':2^{\leq m'}\longrightarrow 2^{<\omega}$ by the following
conditions:
\begin{enumerate}
\item $h'\rest 2^{\leq m(p)}=h$,
\item if $s \in 2^{< m'}\setminus 2^{\leq m(p)}$ then $h'(s)=h(s\rest
m(p))$,
\item if $s \in 2^{m'}$ then $h'(s)=h\lft1(s \rest m(p)\rgt1)^\frown e(s)$.
\end{enumerate}
Thus $m(q)=m'$, but we have to prove that $q=(h',u)$ is a condition. 

Note that for $\<\varepsilon,\phi\>\in u$ we have then:
\begin{multline*}
\sum_{t \in 2^{m'}} 2^{|h'(t)|}\phi\lft1(t,h'(t)\rgt1)= 
\sum_{s\in 2^{m(p)}}\sum_{s \subseteq t\in
2^{m'}} 2^{1+|h(s)|}\cdot \phi(t,h(s)^\frown e(t))>\\
>\sum_{s \in 2^{m(p)}}2^{|h(s)|}\phi\lft1(s,h(s)\rgt1) -
\delta\cdot\sum_{s\in 2^{m(p)}} 2^{1+|h(s)|}>\varepsilon.~\QED
\end{multline*} 

Suppose that $G \subseteq {\cal P}$ is generic over $\V$. Let
$\widetilde{h}_G = \bigcup \{h : \<h,u\> \in G\}$ and  for every $x
\in 2^\omega$, let $h_G(x)=\bigcup_{n \in \omega}
    \widetilde{h}_G(x\rest n)$. It follows
immediately from
\ref{density} that 

  \begin{lemma}
    $h_G(x): 2^\omega \longrightarrow 2^\omega $ is a continuous
    function in $\V[G]$.~$\QED$
  \end{lemma}

\begin{lemma}
  For every measure zero set $H \subseteq 2^\omega\times2^\omega$
  which is coded in $\V$, the set
$$\{x \in 2^\omega: \<x,h_G(x)\> \not \in H\}$$
has measure one. 
\end{lemma}
\Proof
Fix $H$ as above. Suppose that $p=\<h,u\> \in {\cal P}$ and 
$\varepsilon 
>0$ are given. 

It is enough to show that
$\mu\lft2(\{x \in 2^\omega: \<x,h_G(x)\> \not \in H\}\rgt2) >
1-\varepsilon$ holds for every rational $\varepsilon 
>0$. Suppose that $p=\<h,u\> \in {\cal P}$ and $m=m(p)$.

Choose a perfect set $F$ disjoint with $H$ of measure so close to one that
$$\sum_{s\in
2^m} 2^{|h(s)|} \mu([s]\times[h(s)] \cap F) > 1-\varepsilon.$$
Define the function $\phi_F:2^{<\omega}\times2^{<\omega}\rightarrow [0,1]$ by
$$\phi_F(s,t)=\mu\lft1([s]\times [t]\cap F\rgt1)$$
and note that $\<1-\varepsilon, \phi_F\>\in\t$. Moreover, $q=\<h, u
\cup\{\<1-\varepsilon, \phi_F\>\}\>$ is a condition. We show that 
$$q \forces_{\cal P}\mu\lft2(\{x \in 2^\omega: \<x,h_{\dot{G}}(x)\>
\not \in H\}\rgt2) \geq 1-\varepsilon.$$ 
Let $F_n = \bigcup \{[s]\times [t]: s,t \in 2^n \ \& \ 
\lft1([s]\times [t]\rgt1) \cap F \neq \emptyset\}$. Obviously $F =
\bigcap_{n \in \omega} F_n$.
Given $n \in \omega$ there is  $q'=\<h',u'\>$ stronger thatn $q$ such
that $m'=m(q') \geq n$ and $|h'(s)| \geq n $ for $s \in 2^{m'}$. 
Since $\<1- \varepsilon, \phi_F\> \in u \subseteq u'$ also $q''=\<h',
u' \cup \{\<1-\varepsilon, \phi_{F_n}\>\} \>$ is a condition stronger
than $p$ and for $s \in 2^{m'}$, $\lft1([s]\times [h'(s)]\rgt1) \cap
F_n \neq \emptyset$ if and only if $[s]\times [h'(s)] \subseteq F_n$. 
Hence
$$\mu\left( \bigcup \{[s] : s \in 2^{m'} \ \& \ [s] \times [h'(s)]
\subseteq F_n\}\right) = \sum_{s \in 2^{m'}} 2^{|h'(s)|}
\phi_{F_n}\lft1(s, h'(s)\rgt1) > 1-\varepsilon$$
and
$$q'' \forces_{\cal P} \mu\lft2(\{x: \<x,h_{G}(x)\> \in F_n\}\rgt2)
\geq 1 -\varepsilon.$$
 Using density argument and passing to the limit we get 
$$\mu\lft2(\{x: \<x,h_{G}(x)\> \in F\}\rgt2) \geq
1-\varepsilon.~\QED$$

\begin{lemma}\label{nodom}
There exist centered families $\{{\cal P}_i : i \in I\}$, $I$ countable,
such that $\bigcup_{i \in I} {\cal P}_i$ is dense in $ {\cal P}$ and for
every maximal antichain $\{p_n : n \in \omega\}$ in ${\cal P}$ there exists
a natural number $M(i)$ such that for every condition $q \in {\cal P}_i$
there exists $n \leq M(i)$ such that $q $ and $p_n$ are compatible.

In particular, ${\cal P}$ does not add dominating reals.
\end{lemma}
\Proof
For simplicity we will think of the second coordinates of conditions in $\p$
as finite sequences from $\t$.

Let 
$$I=\{(N,k,h,\langle\varepsilon_i: i<N\rangle): k,N\in\omega,
h:2^m\rightarrow2^{<\omega}, \varepsilon_j\in(0,1)\cap\rationals \mbox{ for }
j<N.\}\] 
For $i=(N,k,h,\langle\varepsilon_i: i<N\rangle)\in I$ let 
\begin{multline*}
{\cal P}_i=\lft2\{\<h,\langle\varepsilon_i,\phi_i\rangle_{i<N}\>\in\p: \forall
i<N\ \lft2(\phi_i(\emptyset,\emptyset) \geq 1/k\ \&  \\
  \sum_{s\in
2^m} 2^{|h(s)|} \phi_i\lft1(s,h(s)\rgt1) 
\geq\varepsilon_i+1/k\rgt2)\rgt2\}.
\end{multline*}
Clearly each ${\cal P}_i$ is centered (conditions in $\p$ with the same $h$
can be put together) and they cover $\p$.

We want to show that the families $\p_i$ have the required property.  Assume
not. Thus we have a maximal antichain $\langle p_k:k\in\omega\rangle$ in
$\p$ and a sequence $\langle q_n:n\in\omega\rangle\subseteq\p_i$ (for some
$i=(N,k,h,\langle\varepsilon_i: i<N\rangle)$) such that
$q_n\incomp_{\p} p_k$ for $k  \leq n$.

Let $q_n=\<h,\langle\varepsilon_i, \phi^n_i\rangle_{i<N}\>$,
$\phi^n_i(\emptyset,\emptyset)\geq 1/k$, $\sum_{s \in 2^m} 2^{|h(s)|}
\phi^n_i\lft1(s, h(s)\rgt1)\geq
\varepsilon_i+1/k$. Passing to a subsequence we may assume that for each
$i<N$ the sequence $\langle\phi^n_i: n\in\omega\rangle$ is pointwise
converging (note that the space $[0,1]^\omega$ is compact). 

Let
$\phi_i:2^{<\omega}\times2^{<\omega}\longrightarrow [0,1]$ be the
limit functions, i.e. 
$$\phi_i(s,t)=\lim_{n\to\infty}\phi^n_i(s,t).$$
The functions $\phi_i$ satisfy conditions (1)--(3) of the definition of $\t$
(for the first condition remember that
$\phi^n_i(\emptyset,\emptyset)\geq 1/k$). Moreover
$$\sum_{s\in
2^m} 2^{|h(s)|} \phi_i\lft1(s,h(s)\rgt1)=\lim_{n\to\infty} \sum_{s\in
2^m} 2^{|h(s)|} \phi^n_i\lft1(s,h(s)\rgt1) \geq\varepsilon_i+1/k.$$
Consequently $\<h,\langle\varepsilon_i,\phi_i\rangle_{i<N}\>\in\p$. We find
$k_0\in\omega$ such that the conditions $p_{k_0}$ and
$\<h,\langle\varepsilon_i,\phi_i\rangle_{i<N}\>$ are compatible. Let
$\<h^\star,\langle\varepsilon_i,\phi_i\rangle_{i<N^\star}\>
\geq\<h,\langle\varepsilon_i, 
\phi_i\rangle_{i<N}\>, p_{k_0}$ where $h^\star:2^{\leq m^\star}\longrightarrow
2^{<\omega}$, $N^\star>N$. Then we have for $i<N$:
$$\varepsilon_i< \sum_{s\in
2^{m^\star}} 2^{|h^\star(s)|} \phi_i\lft1(s,h^\star(s)\rgt1) =
\lim_{n\to\infty}  \sum_{s\in
2^{m^\star}} 2^{|h^\star(s)|} \phi^n_i\lft1(s,h^\star(s)\rgt1).$$
Consequently for sufficiently large $n$ we will have
$$\varepsilon_i<\sum_{s\in
2^{m^\star}} 2^{|h^\star(s)|} \phi^n_i\lft1(s,h^\star(s)\rgt1).$$
So take $n>k_0$ such that the above holds for each $i<N$. Then
$\<h^\star,\langle\varepsilon_i,\phi^n_i\rangle_{i <N}\>$ is a condition in $\p$.
Since it is stronger than $q_n$ and compatible with
$\<h^\star,\langle\varepsilon_i,\phi_i\rangle_{i <N^\star}\>$ we conclude that
$q_n\mbox{$\compatible_{\p}$} p_{k_0}$, which contradicts the choice of $q_n$
($n>k_0$ !!!). 

Let $\{i_n: n \in \omega\}$ be an enumeration of $I$ with infinitely
many repetitions. Suppose that $\forces_{\cal P} \dot{f}\in
\omega^\omega$.

Define a function $g \in \V \cap \omega^\omega$ as $g(k) = M(i_k)$,
where $M(i_k)$ is the number obtained by applying the first part of
the lemma to ${\cal P}_{i_k}$ and to the antichain $p_n = \lbv \dot{f}(k)=n\rbv_{\cal P}$, $n
\in \omega$. It is clear that 
$$\forces_{\cal P} \exists^\infty n \ \dot{f}(k) \leq g(k).~\QED$$

\begin{theorem}
  It is consistent with $\ZFCa$ that $\cov^\star(\N)> \max\{\cov(\N),
  \frak b\}$ and that $\unif^\star(\N) < \min\{\unif(\N),
  \frak d\}$.
\end{theorem}
\Proof
To construct the first model let ${\cal P}_{\omega_2}$ be the finite
support iteration of ${\cal P}$ of length $\omega_2$. 

Let $\V \thinks 2^{\boldsymbol\aleph_0}=\boldsymbol\aleph_1$. It is
clear that $\V^{{\cal P}_{\omega_2}} \thinks
\cov^\star(\N)=\boldsymbol\aleph_2$. Since ${\cal P}$ is
$\sigma$-centered neither ${\cal P}$ nor a finite support iteration
of ${\cal P}$ adds random reals (see \cite{shelahbook} or
\cite{BJbook}). Similarly, property stated in \ref{nodom} implies that
finite support iteration of ${\cal P}$ does not add dominating reals.
Thus $\cov(\N)$ and $\frak b$ are both equal to $\boldsymbol\aleph_1$
in $\V^{{\cal P}_{\omega_2}}$. 

The second part of the theorem is proved similarly.  Let $\V \thinks
{\bold { MA}} \ \& \ 2^{\boldsymbol\aleph_0}=\boldsymbol\aleph_2$ and
let ${\cal P}_{\omega_1}$ be the finite support iteration of ${\cal
  P}$ of length $\boldsymbol\aleph_1$.

By ``dualizing'' the above argument we show that
$$\V^{{\cal P}_{\omega_1}} \thinks
\unif^\star(\N)=\boldsymbol\aleph_1 \ \& \ \unif(\N)=\frak
d=\boldsymbol\aleph_2.~\QED$$

\begin{theorem}
Any of the inequalities  $\cov^\star(\N)> {\frak b}$,
$\unif^\star(\N) < {\frak b}$ is consistent with $\ZFCa$.
\end{theorem}
\Proof 
For the first model add $\boldsymbol\aleph_2$ random reals (simultanously) to a model
of $\CH$. Then, in the extension we will have ${\frak d}=\boldsymbol\aleph_1$ and
$\cov^\star(\N)=\boldsymbol\aleph_2$ (for the last note that if $r$ is a random real
over $\V$ then the constant function $h(x)=r$ ``omits'' all measure zero
subsets of the plane coded in $\V$). The second model can be obtained by
adding $\boldsymbol\aleph_1$ random reals to a model of ${\bold {MA}} +
2^{\boldsymbol\aleph_0}=\boldsymbol\aleph_2$.~$\QED $

For the sake of completeness of the picture let us mention the following
result which will appear in \cite{470} (the forcing notion applied for it
is a special case of the scheme presented there):
\begin{theorem}
It is consistent with $\ZFCa$  that  $\cov^\star(\N)<\unif(\M)$.
\end{theorem}

\ifx\undefined\bysame
\newcommand{\bysame}{\leavevmode\hbox to3em{\hrulefill}\,}
\fi

\end{document}